\documentclass[a4paper,10pt]{article}
\usepackage{color}
\begin{document}
\newtheorem{def1}{Definition}[section]
\newtheorem{rem}{Remark}[section]
\newtheorem{prop}{Proposition}[section]
\newtheorem{lem2}{Lemma}[section]
\newtheorem{exam}{Example}[section]
\title{Products of rough finite state machines}
\author{S.P. Tiwari\thanks{sptiwarimaths@gmail.com} and Shambhu Sharan\thanks{shambhupuremaths@gmail.com}\\Department of Applied Mathematics \\Indian School of Mines \\Dhanbad-826004, India}
\date{}
\maketitle
\abstract{In this paper, we introduce the concept of several products of rough finite state machines. We establish their relationships through coverings and investigate some algebraic properties for these products.}\\[5pt]
\section{Introduction} The concept of finite state semiautomata (finite state machines) is well known (cf., e.g., \cite{Gin,Hol,Hop,Ito,Sip}). A (deterministic) finite state machine is  a triple $(Q,X,\delta)$, consisting of two (finite) sets $Q$ (of {\bf states}) and $X$ (of {\bf inputs}) and a map $\delta:Q\times X\rightarrow Q$ (called the {\bf transition map}). A nondeterministic version of a finite state machine, known as nondeterministic finite state machine, is also a triple $(Q,X,\delta)$, where $Q$ and $X$ are as above and $\delta:Q\times X\rightarrow 2^Q$ is a map. The only difference between a deterministic and a nondeterministic finite state machine is in the value that the transition map returns. In case of previous the transition map returns a single state, while in case of later it returns a set of states. \\[5pt]
An account of the fuzzy theoretic version of the notion of a finite state semiautomaton has been studied by Mordeson and Malik in \cite{Mor1}, who called the resulting concept as a fuzzy finite state machine (see also \cite{Mal}). This fuzzy theoretic version is obtained by allowing $\delta(q, a)$, where $q\in Q$ and $a\in X$, to be not just a single state or even a subset of $Q$, but a {\em fuzzy} (sub)set of $Q$. Also, similar, or closely related, notions have been introduced and studied by Kim, Kim and Cho \cite{Kim}, Jun \cite{Jun1}, and Li and Pedrycz \cite{Li}. In literature (c.f., \cite{Cho,Kim,Kum,Mal2,Mor1}), the crisp concepts of several types of products for finite state machines introduced and studied in \cite{Hol} has been fuzzified by many researchers.\\[5pt]
Pawlak's rough set theory \cite{Paw1}, like fuzzy set theory, is another mathematical approach to deal with imprecise,
uncertain or incomplete information and knowledge. It has rapidly drawn attention of both mathematicians and computer scientists due to its ability to model many aspects of artificial intelligence and cognitive sciences, particularly in the areas of knowledge acquisition, decision analysis and expert systems. Following the advent of rough set theory, Basu \cite{Basu} recently introduced the concept of a {\em rough} finite state (semi)automaton, by allowing a state, when given an input, to `transition' to a rough set of the state set (rather than a subset or a fuzzy set) in a certain way and extended the idea further by designing a {\em recognizer} that accepts imprecise statements (cf., \cite{Basu}, for more details). Inspired from the work of Basu, Tiwari and Sharan \cite{Sham} introduced the concepts of rough transformation semigroup associated with a rough finite state machine and coverings of rough finite state machines. Recently, Tiwari, Srivastava and Sharan \cite{Tiw2} introduced and studied the algebraic concepts such as separatedness, connectedness and retrievability of such machines.\\[5pt]
In this present work, our motive is to produce a new rough finite state machine by connecting the rough finite state machines. This we achieve by introducing different types of products between rough finite machines. Specifically, after providing a detail study of rough finite state machines, we introduce and study the notions of different types of products viz., direct product, cascade product and wreath product of rough finite state machines. We explore the relationship among such products through coverings, and investigate some algebraic properties of these products. 
\section{Preliminaries} In this section, we recall and study some concepts associated with rough sets, rough finite state machines and coverings, which we need in the subsequent sections.
\subsection{Rough Sets}
Over the past three decades, a number of definitions of a rough set have appeared in the literature (cf, e.g., \cite{Ban1,Jar,Kom,Pag,Paw1,Paw2,Pol,Yao}). In \cite{Ban}, it has been shown that some of these are equivalent. In this paper, we follow the definition of a rough set as it is given in \cite{Yao}. For completeness, we recall the following key notions.
\begin{def1} {\rm\cite{Paw1}} An \textbf{approximation space} is a pair $(X,R)$, where $X$ is a nonempty set and $R$ is an equivalence relation on $X$.
\end{def1}
If $R$ is an equivalence relation on a nonempty set $X$ and $x\in X$, then let $[x]$ denote the set $\{y\in X\mid {x}R{y}\}$, called an equivalence class or a block under $R$. Let $X/R=\{[x]\mid x\in X\}$.
\begin{def1}\label{def1:l} {\rm\cite{Yao}} Given an approximation space $(X,R)$ and $A \subseteq X$, the \textbf{lower approximation} $\underline{A}$ of $A$ and the \textbf{upper approximation} $\overline{A}$ of $A$ are defined as follows:
\begin{center}
    $\underline{A}=\bigcup\{[x] \in X/R \mid [x]\subseteq A\}$,\\[5pt]
    $\overline{A}=\bigcup\{[x] \in X/R \mid [x]\cap A \neq \phi \}$.
\end{center}
\end{def1}
The pair $(\underline{A},\overline{A})$ is called a \textbf{rough set}. We shall denote it by {\bf A}.\\[5pt]
We now recall another concept of rough set from \cite{Pag}.
\begin{def1}\label{def1:m}  For an approximation space $(X,R)$ and $A \subseteq X$, the pair $(\underline{A},\overline{A}^c)$ is called a \textbf{rough set}.
\end{def1}
Let $(X,R)$ be an approximation space. Define a relation $\equiv$ on $2^X$ by $A\equiv B \Leftrightarrow \overline{A}=\overline{B}$ and $\underline{A}=\underline{B}$. Then $\equiv$ is an equivalence relation on $X$. The following is also a concept of rough set induced by the equivalence relation $\equiv$.
\begin{def1} \label{def1:k} {\rm\cite{Paw1}} A \textbf{rough set} in the approximation space $(X,R)$ is an equivalence class of $P(X)/\equiv$.
\end{def1}
\begin{rem} In \cite{Jar}, the above definition of rough set is given in the generalized setup, precisely, $R$ is a binary relation instead of an equivalence relation on a nonempty set $X$.
\end{rem}
Given an approximation space $(X,R)$ and $A\subseteq X$, $\underline{A}$ and $\overline{A}$ are interpreted as the collection of those objects of the domain $X$ that \textit{definitely} and \textit{possibly} belongs to $A$, respectively. Further, $A$ is called \textbf{definable} (or \textbf{exact}) in $(X,R)$ iff $\underline{A} = \overline{A}$. Equivalently, a definable set is a union of blocks under $R$. For any $A\subseteq X$, $\underline{A}$, $\overline{A}$ and $BnA$ are all definable sets in  $(X,R)$.
\begin{rem} In \cite{Ban}, it has shown that the Definitions {\rm \ref{def1:l}, \ref{def1:m} and \ref{def1:k}} introduced by different researchers at different time, are essentially equivalent to each other for a given approximation space $(X,R)$. Even these equivalent definitions provide different algebras by taking different algebraic operations.
\end{rem}
In our case, we will follow the concept of rough sets given in Definition \ref{def1:l}.
\subsection{Rough finite state machines}
The notion of rough finite state machine has been firstly proposed by Basu \cite{Basu}. In this subsection, our aim is to discuss the concept of a rough finite state machine in details. \\[5pt]
Throughout this section, $X^*$ is the set of all {\em words} on $X$ (i.e., finite strings of elements of $X$, which form a monoid under concatenation of strings) including the empty word (which we shall denote by $e$).\\[5pt]
We begin with the following concept of nondeterministic finite state machine.
\begin{def1} A (nondeterministic) \textbf{finite state machine} is a triple $(Q,X,\delta)$, where $Q$ is a nonempty finite set of states, $X$ is a nonempty
finite set of inputs and a map $\delta:Q\times X\rightarrow 2^Q$, called the {\bf transition map} (or more precisely, $\delta$ is a map such that $\delta(q, a)$, where $q\in Q$ and $a\in X$, is a subset of $Q$).
\end{def1}
The transition map $\delta:Q \times X\rightarrow 2^Q$ can be extended to the map \linebreak $\delta^*:Q \times X^*\rightarrow 2^Q$ such that
\begin{enumerate}
 \item [(i)] $\forall q\in Q$, $\delta(q, e)=\{q\}$, and
\item [(ii)] $\forall q\in Q$, $\forall x\in X^*$ and $\forall a\in X$, $\delta(q,xa)=\bigcup\{\delta(p,a): p\in \delta^*(q,x)\}$.\\
\end{enumerate}
A rough finite state machine is a natural generalization of above nondeterministic finite state machine. The difference is only that in case of a rough finite state machine the transition map returns a rough set of states instead of a set of states, as in the case of nondeterministic finite state machine. This roughness arises due to presence of an equivalence relation on its state-set. Formally, a rough finite state machine can be defined as follows:
\begin{def1} A \textbf{rough finite state machine} (or RFSM) is a 4-tuple $M=(Q,R,X,\delta)$, where $Q$ is a nonempty finite set
(the \textbf{set of states} of $M$), $R$  is an equivalence relation on $Q$, $X$ is a nonempty finite set (the \textbf{set of inputs}) and
$\delta:Q \times X\rightarrow \textbf{A}$, where $\textbf{A}=\{(\underline{A},\overline{A}): A\subseteq Q\}$ is a map (called the \textbf{rough transition map}) such that for each $(q,a)\in Q\times X$, $\delta(q, a)=(\underline{A},\overline{A})$ being a rough set in $(Q,R)$ for some $A\subseteq Q$.\\
\end{def1}
We shall denote  $\underline{A}$ and $\overline{A}$ as $\underline{\delta(q,a)}$ and $\overline{\delta(q,a)}$ respectively. Also, throughout, we will write the set of all rough sets $\{(\underline{A},\overline{A}): A\subseteq Q\}$ in the approximation space $(Q,R)$ just as ${\bf A}$.
\begin{exam} \label{exam:a} Consider a RFSM $(Q,R,X,\delta)$, where $Q=\{q_{1},q_{2},q_{3}, q_{4},\linebreak q_{5}\}$, $R$ is an equivalence relation on $Q$ with $Q/R=\{\{q_{1},q_{2}\},\{q_{3},q_{5}\}, \{q_{4}\}\}$, $X=\{a, b\}$ and the rough transition map $\delta$ is given by the following table:\\[8pt]
\begin{tabular}{|c|l|l|l|l|l|}
        \hline
        $Q$ & $\delta(q,a)$ & $\delta(q,b)$\\[5pt]
        \hline
        $q_{1}$ & $(\{q_{1},q_{2}\},\{q_{1},q_{2}\}\cup \{q_{3},q_{5}\})$ & $(\{q_{4}\},\{q_{3},q_{5}\}\cup \{q_{4}\})$\\[5pt]
        $q_{2}$ & $(\phi,\{q_{3},q_{5}\})$ & $(\{q_{3},q_{5}\},\{q_{1},q_{2}\}\cup\{q_{3},q_{5}\}\cup \{q_{4}\})$\\[5pt]
        $q_{3}$ & $(\{q_{3},q_{5}\}\cup \{q_{4}\},\{q_{1},q_{2}\}\cup $ & $(\{q_{1},q_{2}\},\{q_{1},q_{2}\}\cup \{q_{4}\})$\\[5pt]
        $ $ & $\{q_{3},q_{5}\}\cup \{q_{4}\})$ & $ $\\[5pt]
        $q_{4}$ & $(\{q_{4}\},\{q_{1},q_{2}\}\cup \{q_{4}\})$ & $(\{q_{4}\},\{q_{1},q_{2}\}\cup\{q_{3},q_{5}\}\cup \{q_{4}\})$ \\[5pt]
        $q_{5}$ & $(\{q_{1},q_{2}\}\cup \{q_{4}\},\{q_{1},q_{2}\}\cup $ & $(\{q_{1},q_{2}\},\{q_{1},q_{2}\}\cup \{q_{3},q_{5}\})$\\[5pt]
        $ $ & $\{q_{3},q_{5}\}\cup \{q_{4}\})$ & $ $\\
        \hline
\end{tabular}
\begin{center}
    \rm{Table \ref{exam:a}: State Transition Table}\\
\end{center}
\end{exam}
\begin{rem} Let $(Q,X,\delta)$ be a nondeterministic finite state machine and $R$ be an equivalence relation on $Q$ such that for all $A\subseteq Q, \underline{A}=\overline{A}=C$ (say). Then by identifying $\delta$ with the map $\hat{\delta}:Q \times X\rightarrow \textbf{A}$ given by $\hat{\delta}(q,a)=(C,C)$, $\forall (q,a)\in Q\times X$, we see that every nondeterministic finite state machine can be viewed as a RFSM as defined in {\rm \cite{Ito}}.
\end{rem}
Let $(Q,R,X,\delta)$ be an RFSA and  \textbf{D} be the set of all definable sets generated by $R$ over $Q$. Then transition map $\delta$ of a RFSA $(Q,R,X,\delta)$ can be extended to a map $\delta^*:\textbf{D}\times X \rightarrow \textbf{A}$, as we proceed to explain next.
\begin{def1} Let $(Q,R,X,\delta)$ be an RFSM. Then the \textbf{block transition map}  $\delta^{D}:\textbf{D}\times X\rightarrow \textbf{A}$ is defined as follows: $\forall B\in Q/R$ and $\forall a\in X$,
\begin{eqnarray*}
 \delta^{D}(D,a) &=& \left(\underline{\delta^{D}(D,a)}, \overline{\delta^{D}(D,a)}\right), where\\
 \underline{\delta^{D}(D,a)} &=& \bigcup\{\underline{\delta(q,a)}:q\in B\subseteq D, B\in Q/R\}~ and \\[5pt]
 \overline{\delta^{D}(D,a)}  &=& \bigcup\{\overline{\delta(q,a)}:q\in B\subseteq D, B\in Q/R\}.
\end{eqnarray*}
\end{def1}
\begin{exam} \label{exam:c} Consider the RFSM given in Example \ref{exam:a}. Then the block transitions can be evaluated as:
\begin{eqnarray*}
\delta^D(\{q_{1},q_{2}\}\cup \{q_{3},q_{5}\},a) &=& (\underline{\delta^D(\{q_{1},q_{2}\}\cup \{q_{3},q_{5}\},a)},\\
                                                & &\overline{\delta^D(\{q_{1},q_{2}\}\cup \{q_{3},q_{5}\},a)})\\
                                                &=& (\{q_{1},q_{2}\}\cup\{q_{3},q_{5}\}\cup \{q_{4}\},\\
                                                & & \{q_{1},q_{2}\}\cup\{q_{3},q_{5}\}\cup \{q_{4}\}), since\\
\underline{\delta^D(\{q_{1},q_{2}\}\cup \{q_{3},q_{5}\},a)} &=& \bigcup\{\underline{\delta(q,a)}: q\in B\subseteq \{q_{1},q_{2}\}\cup \{q_{3},q_{5}\}\}\\
                                                            &=& \underline{\delta(q_1,a)}\cup \underline{\delta(q_2,a)}\cup \underline{\delta(q_3,a)}\cup \underline{\delta(q_5,a)}\\
                                        &=& \{q_{1},q_{2}\}\cup\{q_{3},q_{5}\}\cup \{q_{4}\}~ and\\
\overline{\delta^D(\{q_{1},q_{2}\}\cup \{q_{3},q_{5}\},a)} &=& \bigcup\{\overline{\delta(q,a)}: q\in B\subseteq \{q_{1},q_{2}\}\cup \{q_{3},q_{5}\}\}\\
                                                            &=& \overline{\delta(q_1,a)}\cup \overline{\delta(q_2,a)}\cup \overline{\delta(q_3,a)}\cup \overline{\delta(q_5,a)}\\
                                        &=& \{q_{1},q_{2}\}\cup\{q_{3},q_{5}\}\cup \{q_{4}\}.
\end{eqnarray*}
The rest of the block transitions can be computed similarly and are given in the following table.\\[5pt]
\begin{tabular}{|c|l|l|l|l|}
        \hline
        $D$ & $\delta^{D}(D,a)$ & $\delta^{D}(D,b)$\\[5pt]
        \hline
        $\{q_{1},q_{2}\}\cup \{q_{3},q_{5}\}$ & $(\{q_{1},q_{2}\}\cup\{q_{3},q_{5}\}\cup \{q_{4}\},$ & $(\{q_{1},q_{2}\}\cup\{q_{3},q_{5}\}\cup \{q_{4}\},$\\[5pt]
	    & $\{q_{1},q_{2}\}\cup\{q_{3},q_{5}\}\cup \{q_{4}\})$ & $\{q_{1},q_{2}\}\cup\{q_{3},q_{5}\}\cup \{q_{4}\})$\\[5pt]
        $\{q_{1},q_{2}\}\cup \{q_{4}\}$ & $(\{q_{1},q_{2}\}\cup\{q_{4}\},\{q_{1},q_{2}\}\cup$ & $(\{q_{3},q_{5}\}\cup\{q_{4}\},\{q_{1},q_{2}\}\cup$\\[5pt]
	    & $\{q_{3},q_{5}\}\cup \{q_{4}\})$ & $\{q_{3},q_{5}\}\cup \{q_{4}\})$\\[5pt]
        $\{q_{3},q_{5}\}\cup \{q_{4}\}$ & $(\{q_{1},q_{2}\}\cup\{q_{3},q_{5}\}\cup \{q_{4}\},$ & $(\{q_{1},q_{2}\}\cup\{q_{3},q_{5}\}\cup \{q_{4}\},$\\[5pt]
        & $\{q_{1},q_{2}\}\cup\{q_{3},q_{5}\}\cup \{q_{4}\})$ & $\{q_{1},q_{2}\}\cup\{q_{3},q_{5}\}\cup \{q_{4}\})$\\[8pt]
        \hline
\end{tabular}
\begin{center}
    \rm{Table \ref{exam:c}: Block Transition Table}\\
\end{center}
\end{exam}
\begin{def1} \label{def1:c} Let $(Q, R,X,\delta)$ be an RFSM. Define $\delta^*:Q \times X^*\rightarrow \textbf{A}$ as follows:
\begin{description}
\item $(i)$ $\delta^*(q,e)= ([q],[q]), \forall q\in Q$,  and
\item $(ii)$ $\forall q\in Q, \forall x\in X^*$ and $\forall a\in X$, $\delta^*(q,xa)= \left(\underline{\delta^*(q,xa)}, \overline{\delta^*(q,xa)}\right)$, where
$\underline{\delta^*(q,xa)}=\underline{\delta^{D}(\underline{\delta^*(q,x)},a)}$ and
$\overline{\delta^*(q,xa)}=\overline{\delta^{D}(\overline{\delta^*(q,x)},a)}$.
\end{description}
\end{def1}
A block transition map can also be extended, as explained next.
\begin{def1} For an  RFSM $(Q,R,X,\delta)$, the block transition map $\delta^{D}:\textbf{D}\times X\rightarrow \textbf{A}$ can be
extended to a map $\delta^{* D}:\textbf{D}\times X^*\rightarrow \textbf{A}$ as follows: $\forall B\in Q/R$ and $\forall x\in X^*$,
\begin{eqnarray*}
 \delta^{*D}(D,x) &=& \left(\underline{\delta^{*D}(D,x)}, \overline{\delta^{*D}(D,x)}\right), where\\
 \underline{\delta^{*D}(D,x)} &=& \bigcup\{\underline{\delta^*(q,x)}:q\in B\subseteq D, B\in Q/R\}~ and \\
 \overline{\delta^{*D}(D,x)} &=& \bigcup\{\overline{\delta^*(q,x)}:q\in B\subseteq D, B\in Q/R\}.
\end{eqnarray*}
\end{def1}
Following is require to prove the extension of rough transition map.
\begin{def1} Let $(Q, R,X,\delta)$ be an RFSM and $D$ be a definable set generated by $R$ over $Q$. Then
\begin{eqnarray*}
 \delta^{*D}(D,xa) &=& \left(\underline{\delta^{*D}(D,xa)},~ \overline{\delta^{*D}(D,xa)}\right), where\\[5pt]
 \underline{\delta^{*D}(D,xa)} &=& \underline{\delta^{D}(\underline{\delta^{*D}(D,x)},a)}~ and\\[5pt]
 \overline{\delta^{*D}(D,xa)} &=& \overline{\delta^{D}(\overline{\delta^{*D}(D,x)},a)},~ \forall x\in X^* and~ \forall a\in X.
\end{eqnarray*}
\end{def1}
\begin{lem2} Let $(Q,R,X,\delta)$ be an RFSM. Then
\begin{eqnarray*}
            \delta^*(q,xy) &=& \left(\underline{\delta^*(q,xy)}, \overline{\delta^*(q,xy)}\right), where\\[5pt]
\underline{\delta^*(q,xy)} &=& \underline{\delta^{*D}(\underline{\delta^*(q,x)},y)}~ and\\[5pt]
\overline{\delta^*(q,xy)}  &=& \overline{\delta^{*D}(\overline{\delta^*(q,x)},y)},\\
\forall q\in Q ~and~ \forall x,y\in X^*.
\end{eqnarray*}
\end{lem2}
\textbf{Proof:} Let $q\in Q$ and $x,y\in X^*$. We prove the result by induction on $|y|=n$. If $n=1$, let $y=a$. Then from Definition \ref{def1:c}
\begin{eqnarray*}
            \delta^*(q,xa) &=& \left(\underline{\delta^*(q,xa)}, \overline{\delta^*(q,xa)}\right), where\\[5pt]
\underline{\delta^*(q,xa)}  &=& \underline{\delta^{D}(\underline{\delta^*(q,x)},a)}~ and \\[5pt]
\overline{\delta^*(q,xa)}  &=& \overline{\delta^{D}(\overline{\delta^*(q,x)},a)}.
\end{eqnarray*}
Thus the result is true for $n=1$. Now, suppose the result is true for all $x\in X^*$ and $y\in X^*$ such that $|y|=n$. Let $y=ua$, where $|u|=n$. Then
\begin{eqnarray*}
\underline{\delta^*(q,xy)}  &=& \underline{\delta^*(q,xua)}\\
                            &=& \underline{\delta^*(q,za)},~ where ~z=xu\\
                            &=& \underline{\delta^{D}(\underline{\delta^*(q,z)},a)}\\
                            &=& \underline{\delta^{D}({\underline{\delta^{*D}(\underline{\delta^*(q,x)},u)},a})}~~(\textrm{by~ induction}).
\end{eqnarray*}
On the other hand
\begin{eqnarray*}
\underline{\delta^{*D}(\underline{\delta^{*}(q,x)},y)} &=& \underline{\delta^{*D}(\underline{\delta^{*}(q,x)},ua)}\\[5pt]
                                                        &=& \underline{\delta^{*D}(D,ua)},~ where~ D=\underline{\delta^{*}(q,x)}\\
                                                        &=& \underline{\delta^{D}(\underline{\delta^{*D}(D,x)},a)}, ~where~ D=\underline{\delta^{*}(q,x)}\\
                                                       &=& \underline{\delta^{D}({\underline{\delta^{*D}(\underline{\delta^*(q,x)},u)},a})}~.
\end{eqnarray*}
Thus ~~~~~~~$\underline{\delta^*(q,xy)}=\underline{\delta^{*D}(\underline{\delta^*(q,x)},y)}$.\\[5pt]
Similarly, ~~$\overline{\delta^*(q,xy)}=\overline{\delta^{*D}(\overline{\delta^{*}(q,x)},y)}$.\\[5pt] Hence the result is true for $|y|=n+1$.\\\\
Keeping the above in mind, it seems reasonable to accept the following also a definition of an RFSM. By the abuse of notation, we shall write $X$, $\delta$ and $\delta^D$ instead of $X^*$, $\delta^*$ and $\delta^{*D}$ respectively.
\begin{def1} \label{def1:d} A \textbf{rough finite state machine} (or RFSM) is a 4-tuple $M=(Q,R,X,\delta)$, where $Q$ is a nonempty finite set
(the \textbf{set of states} of $M$), $R$  is a given equivalence relation on $Q$, $X$ is a monoid (whose elements are the input symbol) and
$\delta:Q \times X\rightarrow \textbf{A}$ is a map (called the \textbf{rough transition map}) such that
\begin{description}
\item $(i)$ $\delta(q,e)= ([q],[q]), \forall q\in Q$,   and
\item $(ii)$ $\forall q\in Q, \forall x,y\in X$, $\delta(q,xy)= \left(\underline{\delta(q,xy)}, \overline{\delta(q,xy)}\right)$, where\\[5pt]
$\underline{\delta(q,xy)}=\underline{\delta^{D}(\underline{\delta(q,x)},y)}$ and
$\overline{\delta(q,xy)}=\overline{\delta^{D}(\overline{\delta(q,x)},y)}$.\\
\end{description}
\end{def1}
Next, we introduce the concept of homomorphism between two rough finite state machines, which is a natural generalization of the same concept associated with finite state machines. In the case of finite state machines, recall that the \emph{homomorphism} between two nondeterministic finite state machines $(Q,X,\delta)$ and $(R,Y,\mu)$ is a pair of maps $f:Q\rightarrow R$ and $g:X\rightarrow Y$ such that $f(\delta(q,x))\subseteq\mu(f(q),g(x))$.
\begin{def1} \label{def1:a} A \textbf{homomorphism} from an RFSM $M_{1}=(Q_{1}, R_{1},\linebreak X_{1},\delta_{1})$ to a RFSM $M_{2}=(Q_{2},R_{2},X_{2},\delta_{2})$ is a pair of maps $f:Q_{1}\rightarrow Q_{2}$ and $g:X_{1} \rightarrow X_{2}$ such that\\[5pt]
$(i)$ $(p,q) \in R_1 \Rightarrow$ $ (f(p), f(q))\in R_2$, $\forall p,q \in Q_1$, and \\[5pt]
$(ii)$ $f(\delta_1(q,x))\subseteq \delta_2(f(q),g(x))$ or $(f(\underline{\delta_1(q,x)}), f(\overline{\delta_1(q,x)}))\subseteq (\underline{\delta_2(f(q),g(x))}, \linebreak \overline{\delta_2(f(q),g(x))})$, $ \forall q \in Q_1$ and $\forall x\in X_1$.
\end{def1}
A bijective homomorphism $(f,g)$ from an RFSM $M_1$ to an RFSM $M_2$ is called an \textbf{isomorphism}. If there is an isomorphism from RFSM $M_1$ to RFSM $M_2$, then $M_1$ is said to be {\bf isomorphic} to $M_2$, and is denoted by $M_1\cong M_2$.
\begin{exam} \label{exam:f} Let $M_{1}=(Q_{1},R_{1},X_{1},\delta_{1})$ and $M_{2}=(Q_{2},R_{2},X_{2},\delta_{2})$ be two rough finite state machines, where $Q_{1}=\{q_{1},q_{2},q_{3},q_{4}\}$, $Q_{1}/R_{1}=\{\{q_{1}\}, \{q_{2},q_{4}\},\{q_{3}\}\}$, $X_{1}=\{a, b\}$, $Q_{2}=\{q'_{1},q'_{2},q'_{3},q'_{4}\}$, $Q_{2}/R_{2}=\{\{q_{1}'\},\linebreak \{q_{2}',q_{4}'\},\{q_{3}'\}\}$, $X_{2}=\{a', b'\}$ and the rough transition functions $\delta_{1}$ and $\delta_{2}$ are respectively given as follows:\\[5pt]
\begin{tabular}{|c|l|l|l|l|l|}
        \hline
        $Q$ & $\delta_1(q,a)$ & $\delta_1(q,b)$\\[5pt]
        \hline
        $q_{1}$ & $(\{q_{1}\},\{q_{1}\}\cup \{q_{3}\})$ & $(\phi, \{q_{2},q_{4}\})$\\[5pt]
        $q_{2}$ & $(\{q_{3}\},\{q_{2},q_{4}\}\cup \{q_{3}\})$\,\,\,\,\, & $(\{q_{2},q_{4}\},\{q_{1}\}\cup \{q_{2},q_{4}\}\cup \{q_{3}\})$\,\,\,\,\,\,\,\,\,\,\,\,\,\\[5pt]
        $q_{3}$ & $(\phi,\{q_{1}\})$ & $(\{q_{3}\},\{q_{1}\}\cup \{q_{3}\})$\\[5pt]
        $q_{4}$ & $(\{q_{2},q_{4}\},\{q_{1}\}\cup \{q_{2},q_{4}\})$\,\, & $(\{q_{1}\},\{q_{1}\}\cup \{q_{3}\})$ \\[5pt]
        \hline
\end{tabular}
\begin{center}
and
\end{center}
\begin{tabular}{|c|l|l|l|l|l|}
        \hline
        $Q$ & $\delta_2(q',a')$ & $\delta_2(q',b')$\\[5pt]
        \hline
        $q'_{1}$ & $(\{q'_{1}\},\{q'_{1}\}\cup \{q'_{3}\})$ & $(\phi, \{q'_{2},q'_{4}\})$\\[5pt]
        $q'_{2}$ & $(\{q'_{2},q'_{4}\},\{q'_{1}\}\cup \{q'_{2},q'_{4}\})$\,\, & $(\{q'_{1}\},\{q'_{1}\}\cup \{q'_{3}\})$\\[5pt]
        $q'_{3}$ & $(\phi,\{q'_{1}\})$ & $(\{q'_{3}\},\{q'_{1}\}\cup \{q'_{3}\})$\,\,\,\,\,\,\,\,\,\,\,\,\,\,\,\,\,\,\,\,\,\,\\[5pt]
        $q'_{4}$ & $(\{q'_{3}\},\{q'_{2},q'_{4}\}\cup \{q'_{3}\})$ & $(\{q'_{2},q'_{4}\},\{q'_{1}\}\cup \{q'_{2},q'_{4}\}\cup \{q'_{3}\})$\,\,\,\,\,\,\,\,\,\,\,\,\, \\[5pt]
        \hline
\end{tabular}
\begin{center}
    \rm{Table \ref{exam:f}: State Transition Table}\\
\end{center}
A pair of maps $f:Q_{1}\rightarrow Q_{2}$ and $g:X_{1} \rightarrow X_{2}$, where $f(q_{1})=q'_{1}, f(q_{2})=q'_{4}, f(q_{3})=q'_{3}, f(q_{4})=q'_{2}$ and $g(a)=a', g(b)=b'$ is clearly a homomorphism from $M_{1}$ to $M_{2}$.
\end{exam}
\begin{rem} From the Definition \ref{def1:a}, one can easily see that how in a simple way we are introducing the concept of homomorphism in the case of rough finite state machines from the concept of homomorphism of finite state machines. Contrary to it, it is easy to see that if the some concept is known in the case of rough finite state machines, one can easily guess the similar concept in the case of finite state machines. So, now onward we will introduce the concepts for rough finite state machines without recalling the similar concepts for finite state machines.
\end{rem}
The concept of coverings of finite state machines has been introduced and studied in \cite{Hol}. We close this subsection by recalling the concept of covering of  rough finite state machines, recently introduced in \cite{Sham}.
\begin{def1} \label{def1:b} Let $M_1=(Q_1,R_1,X_1,\delta_1)$ and $M_2=(Q_2,R_2,X_2,\delta_2)$ be rough finite state machines. Then a pair of maps  $\eta:Q_2\rightarrow Q_1$ (onto) and $\xi:X_1\rightarrow X_2$ is called a \textbf{covering} of $M_1$ by $M_2$, if
\begin{description}
  \item $(i)$  $(p,q)\in R_2\Rightarrow(\eta(p),\eta(q))\in R_1,\forall p,q \in Q_2$, and
  \item $(ii)$ $\forall q_2\in Q_2$ and $\forall x\in X_1$, $\delta_1(\eta(q_2),x)\subseteq \eta(\delta_2(q_2,\xi(x)))$ or $(\underline{\delta_1(\eta(q_2),x)},\linebreak \overline{\delta_1(\eta(q_2),x)})\subseteq (\eta(\underline{\delta_2(q_2,\xi(x))}),\eta(\overline{\delta_2(q_2,\xi(x))})$,  where $\xi:X_1 \rightarrow X_2$ is a map such that $\xi(e_1)=e_2$ and $\xi(x)=\xi(x_1)\xi(x_2)...\xi(x_n), \forall x=x_1x_2...x_n\in X_1$.
\end{description}
\end{def1}
We shall denote by $M_1\preceq M_2$, the covering of $M_1$ by $M_2$.
\section{Products of rough finite state machines} In this section, we introduce several products for rough finite state machines. We explore the notions of coverings for these products and also examine some algebraic properties. For the terminology in (crisp) automata theory, we refer to \cite{Hol}.\\[5pt]
$R$ appearing below is a relation on $Q_{1}\times Q_{2}$ defined as $((p_1,p_2),(q_1,q_2))\in R$  iff $(p_{1},q_{1})\in R_1$ and $(p_{2},q_{2})\in R_2$. It is easy to see that $R$ turns out to be an equivalence relation on $Q_1 \times Q_2$, as $R_1$ and $R_2$ are equivalence relations on $Q_1$ and $Q_2$ respectively. It is easy to see that the relation $R$ on $Q_{1}\times Q_{2}$ is nothing but $R_1\times R_2$.\\[5pt]
We begin with the following concept of (full) direct product of two rough finite state machines from \cite{Tiw2}. In case of finite state machines, this product may be interpreted as the `parallel composition' of two finite state machines (cf., e.g., D$\ddot{o}$rfler \cite{Dor}).\\
\begin{def1} {\rm\cite{Tiw2}} Let $M_1=(Q_1,R_1,X_1,\delta_1)$ and $M_2=(Q_2,R_2,X_2,\delta_2)$ be rough finite state machines. Then the RFSM $M_1\times M_2=(Q_1\times Q_2,R,X_1\times X_2,\delta_1\times \delta_2)$ is called \textbf{(full) direct product} of $M_1$ and $M_2$, where $\delta_1\times \delta_2:(Q_1\times Q_2)\times (X_1\times X_2)\rightarrow \textbf{A}$, is a map such that $(\delta_1\times \delta_2)((q_1,q_2),(x_1,x_2))=((\underline{\delta_1(q_1,x_1)},\underline{\delta_2(q_2,x_2)}),(\overline{\delta_1(q_1,x_1)},\overline{\delta_2(q_2,x_2)})$, \linebreak $\forall (q_1,q_2)\in Q_1\times Q_2$ and $\forall (x_1,x_2)\in X_1\times X_2$.
\end{def1}
Inspired from \cite{Hol}, we now introduce more direct products of two rough finite state machines.
\begin{def1} Let $M_1=(Q_1,R_1,X,\delta_1)$ and $M_2=(Q_2,R_2,X,\delta_2)$ be rough finite state machines. Then the RFSM $M_1\wedge M_2=(Q_1\times Q_2,R,X,\linebreak \delta_1\wedge \delta_2)$ is called the \textbf{restricted direct product} of $M_1$ and $M_2$, where $\delta_1\wedge \delta_2:(Q_1\times Q_2)\times X\rightarrow \textbf{A}$, is a map such that $(\delta_1\wedge \delta_2)((q_1,q_2),x)=((\underline{\delta_1(q_1,x)},\underline{\delta_2(q_2,x)}), (\overline{\delta_1(q_1,x)},\overline{\delta_2(q_2,x)})$, $\forall (q_1,q_2)\in Q_1\times Q_2$ and $\forall x\in X$.\\
\end{def1}
Let $\overline{X}$ be a finite set and $f:\overline{X}\rightarrow X_1\times X_2$ be a map. Also, let $p_1$ and $p_2$ be the projection mappings of $X_1\times X_2$ onto $X_1$ and $X_2$ respectively, i.e., $p_1:X_1\times X_2\rightarrow X_1$ and $p_2:X_1\times X_2\rightarrow X_2$. Then the following is the concept of generalized direct product of rough finite state machines.
\begin{def1} Let $M_1=(Q_1,R_1,X_1,\delta_1)$ and $M_2=(Q_2,R_2,X_2,\delta_2)$ be rough finite state machines. Then the RFSM $M_1\ast M_2=(Q_1\times Q_2,R,\overline{X},\delta_1\ast \delta_2)$ is called \textbf{general direct product} of $M_1$ and $M_2$, where $\delta_1\ast \delta_2:(Q_1\times Q_2)\times \overline{X}\rightarrow \textbf{A}$, is a map such that $(\delta_1\ast \delta_2)((q_1,q_2),x)=((\underline{\delta_1(q_1,p_1(f(x)))}, \underline{\delta_2(q_2,p_2(f(x)))}),(\overline{\delta_1(q_1,p_1(f(x)))},
\overline{\delta_2(q_2,p_2(f(x)))}))$, \linebreak $\forall (q_1,q_2)\in Q_1\times Q_2$ and $\forall x\in \overline{X}$.
\end{def1}
\begin{rem} $(i)$ If $\overline{X}=X_1\times X_2$ and $f$ is the identity map, then the general direct product $M_1\ast M_2$ reduces to full direct product. \\[5pt]
$(ii)$ If $\overline{X}=X_1=X_2$ and $f$ is the identity map, then the general direct product $M_1\ast M_2$ reduces to restricted direct product.
\end{rem}
The following proposition shows the relation between (full) direct product and restricted direct product through covering.
\begin{prop} Let $M_1=(Q_1,R_1,X,\delta_1)$ and $M_2=(Q_2,R_2,X,\delta_2)$ be rough finite state machines. Then $M_1\wedge M_2\preceq M_1\times M_2$.
\end{prop}
\textbf{Proof.} Let $\eta$ be an identity map on $Q_1\times Q_2$. Then $\forall ((q_1,q_2), (q_1',q_2'))\in Q_1\times Q_2$, $((q_1,q_2),(q_1',q_2'))\in R \Rightarrow (\eta(q_1,q_2),\eta(q_1',q_2'))\in R$. Define a map $\xi:X\rightarrow X\times X$ by $\xi(x)=(x,x)$, $\forall x\in X$. Now, for $(q_1,q_2)\in Q_1\times Q_2$ and $x\in X$, $(\delta_1\wedge \delta_2)(\eta(q_1,q_2),x)=(\delta_1\wedge \delta_2)((q_1,q_2),x)= ((\underline{\delta_1(q_1,x)},\underline{\delta_2(q_2,x)}),  (\overline{\delta_1(q_1,x)}, \overline{\delta_2(q_2,x)})=(\delta_1\times \delta_2)((q_1,q_2),(x,x))=(\delta_1\times \delta_2)((q_1,q_2),\xi(x))$. Also, $\xi(e)=(e,e)$, $e$ being the identity of $X$ and $\xi(x)=\xi(x_1)\xi(x_2).....\xi(x_{n-1}) \xi(x_n), \forall x=x_1x_2...x_n\in X$. Thus $M_1\wedge M_2\preceq M_1\times M_2$.\\\\
The following lemma is useful to introduce the wreath product of rough finite state machines.
\begin{lem2} Let $S_1$ and $S_2$ be semigroups. Then $(S_{1}^{Q_{2}}\times S_2, \ast)$ is a semigroup, where $S_1^{Q_2}=\{f\mid f:Q_2\rightarrow S_1\}$,  $(f, s)\ast (g, t)=(fg, st)$ and  $(fg)(q_2)=f(q_2)g(q_2)$, $\forall f,g\in S_{1}^{Q_{2}}$, $s,t \in S_2$ and $q_2\in Q_2$.
\end{lem2}
\textbf{Proof.}  Let $f,g,h\in S_{1}^{Q_{2}}$ and $s,t,u\in Q_2$. Then $((f,s)\ast (g,t))\ast (h,u)=(fg,st)\ast (h,u)=((fg)h,(st)u)=(f(gh),s(tu))=(f,s)\ast (gh,tu)=(f,s)\ast ((g,t)\ast (h,u))$. Now, let $I\in S_{1}^{Q_{2}}$ such that $I(q_2)=e_1$, $e_1$ being the identity of $S_1$. Then it can be seen that $(I,e_2)$ is an identity of $S_{1}^{Q_{2}}\times S_2$, where $e_2$ is the identity of $S_2$. Thus $(S_{1}^{Q_{2}}\times S_2, \ast)$ is a semigroup with identity $(I,e_2)$.\\[5pt]
Now, we introduce the wreath product of rough finite state machines, which is a generalization of the same concept for finite state machines (cf., \cite{Hol}).
\begin{def1} Let $M_1=(Q_1,R_1,X_1,\delta_1)$, $M_2=(Q_2,R_2,X_2,\delta_2)$ be rough finite state machines. Then the RFSM $M_1 \circ M_2=(Q_1\times Q_2,R,\linebreak X_1^{Q_2}\times X_2,\delta_1 \circ \delta_2)$ is called the \textbf{wreath product} of $M_1$ and $M_2$, where $\delta_1 \circ \delta_2:(Q_1\times Q_2)\times (X_1^{Q_2}\times X_2)\rightarrow \textbf{A}$, is a map such that $(\delta_1\circ\delta_2)((q_1,q_2),(f,x))=((\underline{\delta_1(q_1,f(q_2))},
\underline{\delta_2(q_2,x)}),(\overline{\delta_1(q_1,f(q_2))},\overline{\delta_2(q_2,x)}))$, \linebreak $\forall (q_1,q_2)\in Q_1\times Q_2$ and $\forall (f,x)\in X_1^{Q_2}\times X_2$.\\
\end{def1}
Let $M_n=(Q_n,R_n,X_n,\delta_n), n=1, 2, 3,4$ be rough finite state machines. Then $\delta_i\times \delta_j$ and $\delta_i\circ \delta_j$, $i=1, 2,3,4$, $j=1, 2,3,4,$ appearing below associated with rough finite state machines $M_i\times M_j$ and $M_i\circ M_j$, $i=1,2,3,4$, $j=1,2,3,4$, respectively have their usual meaning.\\[5pt]
Also, $R_i\times R_j$ appearing below is a relation on $Q_i\times Q_j$ defined as $((p_i,p_j),\linebreak (q_i,q_j)) \in R_i\times R_j$ iff $(p_{i},q_{i})\in R_i$ and  $(p_{j},q_{j})\in R_j$, $i,j=1,2,3,4$.\\[5pt]
Now, we have the following.
\begin{prop} Let $M_i=(Q_i,R_i,X_i,\delta_i)$, where $i=1,2,3,4$ be rough finite state machines. Then $(M_1 o M_2)\times (M_3 o M_4)\preceq (M_1\times M_3) o (M_2\times M_4)$.
\end{prop}
{\bf Proof:} Let $M_1 o M_2=(Q_1\times Q_2,R_1\times R_2,(X_1^{Q_2}\times X_2),\delta_1 o \delta_2)$ and $M_3 o M_4=(Q_3\times Q_4,R_3\times R_4,(X_3^{Q_4}\times X_4),\delta_3 o \delta_4)$. Then $(M_1 o M_2)\times (M_3 o M_4)=((Q_1\times Q_2)\times (Q_3\times Q_4), (R_1\times R_2)\times (R_3\times R_4), (X_1^{Q_2}\times X_2)\times (X_3^{Q_4}\times X_4), (\delta_1 o \delta_2)\times (\delta_3 o \delta_4))$. Again, let $M_1\times M_3=(Q_1\times Q_3,R_1\times R_3,X_1\times X_3,\delta_1\times \delta_3)$ and $M_2\times M_4=(Q_2\times Q_4,R_2\times R_4,X_2\times X_4,\delta_2\times \delta_4)$, then $(M_1\times M_3)o (M_2\times  M_4)=((Q_1\times Q_3)\times (Q_2\times Q_4), (R_1\times R_3)\times (R_2\times R_4), (X_1\times X_3)^{{Q_2}\times {Q_4}}\times (X_2\times X_4), (\delta_1\times \delta_3) o (\delta_2 \times \delta_4))$. Define $\eta: (Q_1\times Q_3)\times (Q_2\times Q_4)\rightarrow (Q_1\times Q_2)\times (Q_3\times Q_4)$ by $\eta ((q_1,q_3),(q_2,q_4))=((q_1,q_2),(q_3,q_4))$, for $q_i\in Q_i$, where $i=1,2,3,4$. Then $\eta$ is an onto mapping. Again, define $f\times g: Q_2\times Q_4\rightarrow X_1\times X_3$ by $(f\times g)(q_2,q_4)=(f(q_2),g(q_4))$, where $(q_2,q_4)\in Q_2\times Q_4$, $f:Q_2\rightarrow X_1$ and $g: Q_4\rightarrow X_3$ are functions. Now, define a map $\xi:(X_1^{Q_2}\times X_2)\times (X_3^{Q_4}\times X_4)\rightarrow (X_1\times X_3)^{{Q_2}\times {Q_4}}\times (X_2\times X_4)$ by $\xi((f,x_2),(g,x_4))=(f\times g, (x_2,x_4))$, where $((f,x_2),(g,x_4))\in (X_1^{Q_2}\times X_2)\times (X_3^{Q_4}\times X_4)$. Then $\forall q_i\in Q_i$, where $i=1,2,3,4$, $((\delta_1 o \delta_2)\times (\delta_3 o \delta_4))(\eta((q_1,q_3),(q_2,q_4)),((f,x_2),(g,x_4)))=$ $((\delta_1 o \delta_2)((q_1,q_2), (f(q_2),x_2)),(\delta_3 o \delta_4)((q_3,q_4),(g(q_4),x_4)))=(((\underline{\delta_1(q_1, f(q_2))},\linebreak \underline{\delta_2(q_2,x_2)}),(\overline{\delta_1(q_1,f(q_2))},
\overline{\delta_2(q_2,x_2)}))$, $((\underline{\delta_3(q_3,g(q_4))},\underline{\delta_4(q_4,x_4)}),(\overline{\delta_3(q_3,g(q_4))}$, $
\overline{\delta_4(q_4,x_4)})))= (((\underline{\delta_1(q_1,f(q_2))}, \underline{\delta_3(q_3,g(q_4))}),(\overline{\delta_1(q_1,f(q_2))},
\overline{\delta_3(q_3,g(q_4))})),\linebreak((\underline{\delta_2(q_2,x_2)},  \underline{\delta_4(q_4,x_4)}),(\overline{\delta_2(q_2,x_2)},
\overline{\delta_4(q_4,x_4)})))=((\delta_1\times \delta_3)((q_1,q_3), (f(q_2),\linebreak g(q_4)),(\delta_2\times \delta_4)((q_2,q_4), (x_2,x_4))))=((\delta_1\times \delta_3)o (\delta_2\times \delta_4))(((q_1,q_3),(q_2,q_4),\linebreak ((f\times g),(x_2,x_4))))=((\delta_1\times \delta_3)o (\delta_2\times \delta_4))(((q_1,q_3),(q_2,q_4),\xi((f,x_2),\linebreak(g,x_4))))$.\\[5pt]
Finally, we introduce the following concept of cascade product of rough finite state machines.
\begin{def1} Let $M_1=(Q_1,R_1,X_1,\delta_1)$, $M_2=(Q_2,R_2,X_2,\delta_2)$ be rough finite state machines and $\omega:Q_2\times X_2 \rightarrow X_1$ be a map. Then the RFSM $M_1\omega M_2=(Q_1\times Q_2,R,X_2,\delta_1\omega \delta_2)$ is called the \textbf{cascade product} of $M_1$ and $M_2$, where $\delta_1\omega \delta_2:(Q_1\times Q_2)\times X_2\rightarrow \textbf{A}$, is a map such that $(\delta_1\omega \delta_2)((q_1,q_2),x_2) =((\underline{\delta_1(q_1,\omega(q_2,x_2))},\underline{\delta_2(q_2,x_2)}), (\overline{\delta_1(q_1,\omega(q_2,x_2))}, \linebreak \overline{\delta_2(q_2,x_2)}))$, $\forall (q_1,q_2)\in Q_1\times Q_2$ and $\forall x_2\in X_2$.
\end{def1}
\begin{rem} Let $M_1\omega M_2=(Q_1\times Q_2,R,X_2,\delta_1\omega \delta_2)$ be the cascade product of rough finite state machines $M_1$ and $M_2$ such that $X_1=X_2=X$(say) and $\omega: Q_2\times X\rightarrow X$ be the map, then the restricted direct products of $M_1$ and $M_2$ is a special case of their cascade products.
\end{rem}
Now, we have the following interesting covering property between wreath product and cascade product of rough finite state machines.
\begin{prop} Let $M_1$ and $M_2$ be rough finite state machines. Then $M_1\omega M_2\preceq M_1 \circ M_2$, where $\omega:Q_2\times X_2 \rightarrow X_1$ is a map.
\end{prop}
\textbf{Proof.}  Let $M_1=(Q_1,R_1,X_1,\delta_1)$, $M_2=(Q_2,R_2,X_2,\delta_2)$ be rough finite state machines and $\omega:Q_2\times X_2 \rightarrow X_1$ be a map. Then $M_1\omega M_2=(Q_1\times Q_2,R,X_2,\delta_1\omega \delta_2)$ and $M_1 o M_2=(Q_1\times Q_2,R,X_1^{Q_2}\times X_2,\delta_1 o \delta_2)$, where $\delta_1 \circ \delta_2:(Q_1\times Q_2)\times (X_1^{Q_2}\times X_2)\rightarrow \textbf{A}$, is a map such that $(\delta_1 \circ \delta_2)((q_1,q_2),(f,x))=((\underline{\delta_1(q_1,f(q_2))},\underline{\delta_2(q_2,x)}),  (\overline{\delta_1(q_1,f(q_2))},\overline{\delta_2(q_2,x)}))$, \linebreak $\forall (q_1,q_2)\in Q_1\times Q_2$ and $\forall (f,x)\in X_1^{Q_2}\times X_2$. Let $\eta$ be an identity map on $Q_1\times Q_2$. Then $\forall ((q_1,q_2), (q_1',q_2'))\in Q_1\times Q_2$, $((q_1,q_2),(q_1',q_2'))\in R \Rightarrow (\eta(q_1,q_2),\eta(q_1',q_2'))\in R$. Define a map $\xi:X_2\rightarrow X_1^{Q_2}\times X_2$ by $\xi(x_2)=(f,x_2)$, where  $f:Q_2\rightarrow X_1$ such that $f(q_2)=\omega(q_2,x_2)$, $\forall q_2\in Q_2$. Then for $(q_1,q_2)\in Q_1\times Q_2$ and $x_2\in X_2$, $(\delta_1\omega \delta_2)(\eta(q_1,q_2),x_2)=(\delta_1\omega \delta_2)((q_1,q_2), x_2)=((\underline{\delta_1(q_1,\omega(q_2,x_2))}, \underline{\delta_2(q_2,x_2)}), (\overline{\delta_1(q_1,\omega(q_2,x_2))},  \linebreak \overline{\delta_2(q_2,x_2)}))= ((\underline{\delta_1(q_1,f(q_2))}, \underline{\delta_2(q_2,x_2)}),(\overline{\delta_1(q_1,f(q_2))},  \overline{\delta_2(q_2,x_2)}))=\linebreak (\delta_1 o \delta_2)((q_1,q_2), (f,x_2))= (\delta_1 o \delta_2) ((q_1,q_2), \xi(x_2))$. Again, let $I\in X_1^{Q_2}$ such that $I(q_2)=e_1$, $e_1$ being the identity of $X_1$. Then $\xi(e_2)=(I(q_2),e_2)\linebreak =(e_1,e_2)$ and $\xi(x)=(f,x)=(f_1f_2...f_n,x_1x_2...x_n)=((f_1,x_1)(f_2,x_2)...\linebreak (f_n,x_n))$, $\forall f=f_1f_2...f_n\in X_1^{Q_2}$ and $\forall x=x_1x_2...x_n\in X_2$. Thus $M_1\omega M_2\preceq M_1 \circ M_2$.\\[5pt]
The following propositions are direct consequences of the associativity of products of rough finite state machines.
\begin{prop} Let $M_1, M_2$ and $M_3$ be rough finite state machines. Then
\begin{enumerate}
  \item [(i)] $(M_1\times M_2)\times M_3\cong M_1\times (M_2\times M_3)$,
  \item [(ii)] $(M_1\wedge M_2)\wedge M_3\cong M_1\wedge (M_2\wedge M_3)$,
  \item [(iii)] $(M_1\circ M_2)\circ M_3\cong M_1\circ (M_2\circ M_3)$, and
  \item [(iv)] $(M_1\omega_1 M_2)\omega_2 M_3\cong M_1\omega_3 (M_2\omega_4 M_3)$, where $\omega_3$ and $\omega_4$ are determined by $\omega_1$ and $\omega_2$ in a natural way.\\
\end{enumerate}
\end{prop}
Let $M_n=(Q_n,R_n,X_n,\delta_n), n=1, 2, 3$ be rough finite state machines. Then $\delta_i\times \delta_j$, $\delta_i\wedge \delta_j$, $\delta_i\circ \delta_j$ and $\delta_i\omega\delta_j$, where $\omega_i:Q_j\times X_j\rightarrow X_i$, $i=1,2$, $j=3$ appearing below are rough transition maps associated with rough finite state machines $M_i\times M_j$, $M_i\wedge M_j$, $M_i\circ M_j$ and $M_i\omega_i M_j$, $i=1,2$, $j=3$, respectively.
\begin{prop} \label{prop:k} Let $M_1=(Q_1,R_1,X_1,\delta_1)$, $M_2=(Q_2,R_2,X_2,\delta_2)$ and $M_3=(Q_3,R_3,X_3,\delta_3)$ be rough finite state machines such that $M_1\preceq M_2$. Then
\begin{enumerate}
  \item [(i)] $(a)$ $M_1\times M_3\preceq M_2\times M_3$ and $(b)$ $M_3\times M_1\preceq M_3\times M_2$,
  \item [(ii)] if $X=X_1=X_2$, then $(a)$ $M_1\wedge M_3\preceq M_2\wedge M_3$ and $(b)$ $M_3\wedge M_1\preceq M_3\wedge M_2$,
  \item [(iii)] $(a)$ $M_1\circ M_3\preceq M_2\circ M_3$ and $(b)$ $M_3\circ M_1\preceq M_3\circ M_2$,
  \item [(iv)] $(a)$ given $\omega_1:Q_3\times X_3\rightarrow X_1$ there exists $\omega_2:Q_3\times X_3\rightarrow X_2$, such that $M_1\omega_1 M_3\preceq M_2\omega_2 M_3$ and $(b)$ if $(\eta,\xi)$ is a covering of $M_1$ by $M_2$, then for each $\omega_1:Q_1\times X_1\rightarrow X_3$ there exists $\omega_2:Q_2\times X_2\rightarrow X_3$ such that $M_3\omega_1 M_1\preceq M_3\omega_2 M_2$.
\end{enumerate}
\end{prop}
\textbf{Proof.}  As $M_1\preceq M_2$, there exist an onto map $\eta:Q_2\rightarrow Q_1$ and a map $\xi:X_1\rightarrow X_2$ such that $(q_2,q_2')\in R_2\Rightarrow(\eta(q_2),\eta(q_2'))\in R_1,\forall q_2,q_2' \in Q_2$, and $\delta_1(\eta(q_2),x)\subseteq \eta(\delta_2(q_2,\xi(x)))$ or $ (\underline{\delta_1(\eta(q_2),x)},\overline{\delta_1(\eta(q_2),x)})\subseteq \linebreak (\eta(\underline{\delta_2(q_2,\xi(x))}), \eta(\overline{\delta_2(q_2,\xi(x))}))$, $\forall q_2\in Q_2$ and $\forall x\in X_1$, where $\xi(e)=e$ and $\xi(x)=\xi(x_1)\xi(x_2)...\xi(x_n), \forall x=x_1x_2...x_n \in X_1$.
\\[5pt] $(i)$ $(a)$ Let $M_1\times M_3=(Q_1\times Q_3,R_1\times R_3,X_1\times X_3,\delta_1\times \delta_3)$ and $M_2\times M_3=(Q_2\times Q_3,R_2\times R_3,X_2\times X_3,\delta_2\times \delta_3)$. Define an onto map $\eta_{\times}:Q_2\times Q_3 \rightarrow Q_1\times Q_3$ by $\eta_{\times}(q_2,q_3)=(\eta(q_2),q_3)$ and a map  $\xi_{\times}:X_1\times X_3 \rightarrow X_2\times X_3$ by $\xi_{\times}(x_1,x_3)=(\xi(x_1),x_3)$. Then $((q_2,q_3),(q_2',q_3'))\in R_2\times R_3\Rightarrow ((\eta(q_2),q_3),(\eta(q_2'),q_3'))\in  R_1\times R_3$, $\forall ((q_2,q_3), (q_2',q_3'))\in Q_2\times Q_3$. Let $(q_2,q_3)\in Q_2\times Q_3$ and $(x_1,x_3)\in X_1\times X_3$, then $(\delta_1\times \delta_3)(\eta_{\times}(q_2,q_3),(x_1,x_3))\linebreak =(\delta_1\times \delta_3)(\eta(q_2),q_3),(x_1,x_3))=((\underline{\delta_1(\eta(q_2),x_1)},
\underline{\delta_3(q_3,x_3)}),(\overline{\delta_1(\eta(q_2),x_1)}, \linebreak \overline{\delta_3(q_3,x_3)})) \subseteq
((\eta(\underline{\delta_2(q_2,\xi(x_1))}), {\delta_3(q_3,x_3)}), (\eta(\overline{\delta_2(q_2,\xi(x_1))}), \overline{\delta_3(q_3,x_3)}))=\eta((\delta_2\times \delta_3)(q_2,q_3), (\xi(x_1),x_3))=\eta((\delta_2\times \delta_3)((q_2,q_3), \xi_{\times}(x_1,x_3)))$. Now, $\xi_{\times}(e_1,e_3)=(\xi(e_1),e_3)=(e_1,e_3)$,  where  $e_1$, $e_3$ being the identity elements of $X_1$, $X_3$ respectively and $\xi_{\times}(x,y) =(\xi(x),y) =(\xi(x_1x_2...x_n), y_1y_2...y_n)\linebreak  =(\xi(x_1)\xi(x_2)...\xi(x_n), y_1y_2...y_n)=((\xi(x_1),y_1) (\xi(x_2),y_2)...(\xi(x_n), y_n))$, \linebreak  $\forall x=x_1x_2...x_n\in X_1$ and $\forall y=y_1y_2...y_n \in  X_3$. Thus $M_1\times M_3\preceq M_2\times M_3$.
\\[5pt] $(b)$ The proof is similar to that of Proposition \ref{prop:k} $(i)$ $(a)$.
\\[5pt]$(ii)$ $(a)$ Let $X=X_1=X_2$. Then $M_1\wedge M_3=(Q_1\times Q_3,R_1\times R_3,X, \delta_1\wedge \delta_3)$ and $M_2\wedge M_3=(Q_2\times Q_3,R_2\times R_3,X, \delta_2\wedge \delta_3)$. Define an onto map $\eta_{\wedge}:Q_2\times Q_3 \rightarrow Q_1\times Q_3$ by $\eta_{\wedge}(q_2,q_3)=(\eta(q_2),q_3)$ and $\xi_{\wedge}=\xi$ be the map on $X$. Then $(\eta_{\wedge},\xi_{\wedge})$ is the required covering.
\\[5pt] $(b)$ Follows as above.
\\[5pt] $(iii)$ $(a)$ Let $M_1 \circ M_3=(Q_1\times Q_3,R_1\times R_3,X_1^{Q_3}\times X_3,\delta_1 \circ \delta_3)$ and $M_2 \circ M_3=(Q_2\times Q_3,R_2\times R_3,X_2^{Q_3}\times X_3,\delta_2 \circ \delta_3)$. Define an onto map $\eta_{\circ}:Q_2\times Q_3 \rightarrow Q_1\times Q_3$ by $\eta_{\circ}(q_2,q_3)=(\eta(q_2),q_3)$ and a map $\xi_{\circ}:X_1^{Q_3}\times X_3 \rightarrow X_2^{Q_3}\times X_3$ by $\xi_{\circ}(f,x_3)=(\xi\circ f,x_3)$. Then $((q_2,q_3),(q_2',q_3'))\in R_2\times R_3\Rightarrow ((\eta(q_2),q_3),(\eta(q_2'),q_3'))\in  R_1\times R_3$, $\forall ((q_2,q_3), (q_2',q_3'))\in Q_2\times Q_3$. Also, let $(q_2,q_3)\in Q_2\times Q_3$ and $(f,x_3)\in X_1^{Q_3}\times X_3$. Then $(\delta_1\circ \delta_3)(\eta_{\circ}(q_2,q_3), (f,x_3))= (\delta_1\circ \delta_3)(\eta(q_2),q_3),(f,x_3))=((\underline{\delta_1(\eta(q_2),f(q_2))},\linebreak
\underline{\delta_3(q_3,x_3)}), (\overline{\delta_1(\eta(q_2),f(q_2))},  \overline{\delta_3(q_3,x_3)})) \subseteq
((\eta(\underline{\delta_2(q_2,\xi(f(q_2)))}), \underline{\delta_3(q_3,x_3)}),\linebreak(\eta(\overline{\delta_2(q_2,\xi(f(q_2)))}),
 \overline{\delta_3(q_3,x_3)}))=\eta((\delta_2\circ \delta_3)(q_2,q_3),(\xi\circ f,x_3))=
\eta((\delta_1\circ \delta_3)((q_2,q_3),\xi_{\circ}(f,x_3)))$. Now, for $I\in X_1^{Q_3}$, $I(q_3)=e_1$, $e_1$ being the identity of $X_1$, $(I, e_3)$ is an identity element of $X_1^{Q_3}\times X_3$. Also, $\xi_{\circ}(I,e_3)=(\xi\circ I,e_3)$, which is an identity element of $X_2^{Q_3}\times X_3$. Again, $\xi_{\circ}(f,x)=(\xi\circ f,x)=(\xi\circ (f_1f_2...f_n), x_1x_2...x_n)=((\xi\circ f_1)(\xi\circ f_2)...(\xi\circ f_n),x_1x_2...x_n)=((\xi\circ f_1,x_1)(\xi\circ f_2,x_2)...(\xi\circ f_n,x_n)$, $\forall f=f_1f_2...f_n\in X_1^{Q_3}$ and $x=x_1x_2...x_n\in X_3$. Thus $M_1\circ M_3\preceq M_2\circ M_3$.
\\[5pt] $(b)$ The proof is similar to that of Proposition \ref{prop:k} $(iii)$ $(a)$.
\\[5pt]$(iv)$ $(a)$ For given $\omega_1:Q_3\times X_3\rightarrow X_1$, let $M_1\omega_1 M_3=(Q_1\times Q_3,R_1\times R_3,X_3,\delta_1\omega_1 \delta_3)$. Then  there exists $\omega_2:Q_3\times X_3\rightarrow X_2$ such that $M_2\omega_2 M_3=(Q_2\times Q_3,R_2\times R_3,X_3,\delta_2\omega_2 \delta_3)$. Define an onto map $\eta_{\omega}:Q_2\times Q_3 \rightarrow Q_1\times Q_3$ such that $\eta_{\omega}(q_2,q_3)=(\eta(q_2),q_3)$, $\forall(q_2,q_3)\in Q_2\times Q_3$. Then $((q_2,q_3),(q_2',q_3'))\in R_2\times R_3\Rightarrow ((\eta(q_2),q_3),(\eta(q_2'),q_3'))\in  R_1\times R_3$, \linebreak $\forall ((q_2,q_3), (q_2',q_3'))\in Q_2\times Q_3$. Now, let  $\omega_2=\omega_1$ and $\xi_{\omega}$ be an identity map on $X_3$. Also, let $(q_2,q_3)\in Q_2\times Q_3$ and $x_3\in X_3$. Then $(\delta_1\omega_1 \delta_3)(\eta_{\omega}(q_2,q_3),x_3) =(\delta_1\omega_1\delta_3)((q_2,q_3),x_3)=\linebreak
((\underline{\delta_1(\eta(q_2), \omega_1(q_3,x_3))}, \underline{\delta_3(q_3,x_3)}), (\overline{\delta_1(\eta(q_2),\omega_1(q_3,x_3))}, \overline{\delta_3(q_3,x_3)}))\subseteq \linebreak((\eta(\underline{\delta_2(q_2,\omega_2(q_3,x_3))}),\underline{\delta_3(q_3,x_3)}),
(\eta(\overline{\delta_2(q_2,\omega_2(q_3,x_3))}),\overline{\delta_3(q_3,x_3)}))=\linebreak \eta((\delta_2\omega_2\delta_3)((q_2,q_3),x_3)$. Hence the covering exist.
\\[5pt] $(b)$ Given $\omega_1:Q_1\times X_1\rightarrow X_3$, let $\omega_2:Q_2\times X_2\rightarrow X_3$ such that $\omega_2(q_2,\xi(x_1))=\omega_1(\eta(q_2),x_1)$. Define an onto map $\eta_{\omega}:Q_3\times Q_2 \rightarrow Q_3\times Q_1$ by $\eta_{\omega}(q_3,q_2)=(q_3,\eta(q_2))$ and set $\xi_{\omega}=\xi$. Then $((q_3,q_2),(q_3',q_2')) \in R_3\times R_2\Rightarrow ((q_3,\eta(q_2)),((q_3',\eta(q_2'))\in  R_3\times R_1$, $\forall ((q_3,q_2), (q_3',q_2'))\in Q_3\times Q_2$. Thus $(\eta_{\omega}, \xi{\omega})$ is the required covering.
\section{Conclusion} Chiefly inspired from \cite{Hol} and \cite{Mor1}, we have introduced and studied here the concept of rough finite state machine and several products viz., direct product, cascade product and wreath product of rough finite state machines. Also, we studied the relationship between these different products through coverings as well as examined some algebraic properties. We hope that, like fuzzy finite state machines, rough finite state machines, which is another dimension of application of rough set theory, will attract the researchers and the work carried out here will help in finding some successful applications of rough finite state machines.

\end{document}